\documentclass[11pt]{article}
\usepackage{amsmath}
\usepackage{bbm}
\usepackage{amsmath,amsthm,amssymb,amscd}
\usepackage{latexsym}
\usepackage{hyperref}
\usepackage{indentfirst}
\usepackage[numbers,sort&compress]{natbib}

\newtheorem{theorem}{Theorem}[section]
\newtheorem{corollary}{Corollary}[section]
\newtheorem{proposition}{Proposition}[section]
\newtheorem{lemma}{Lemma}[section]
\newtheorem{definition}{Definition}[section]
\newtheorem{remark}{Remark}[section]

\numberwithin{equation}{section} 

\title{\bf Liouville type theorems for stable solutions of the weighted elliptic system}
\author{{\small Liang-Gen Hu$^1$\hspace*{10pt}Jing Zeng$^2$}\\[0.05cm] {\small 1. Department of Mathematics,
Ningbo University, 315211, Ningbo, P.R. China}\\ {\small email address: hulianggen@tom.com}\\ {\small 2. School of Mathematics and Computer Science,
Fujian Normal University,}\\ {\small 350007, Fuzhou,  P.R. China, email address: zengjing@fjnu.edu.cn}}
\date{}
\begin{document}
\maketitle
\def\abstractname{}
\vspace*{-40pt}

\begin{abstract}
\noindent {\bf Abstract:} We examine the weighted elliptic system
\begin{equation*}
\begin{cases}
-\Delta u=(1+|x|^2)^{\frac{\alpha}{2}} v,\\
-\Delta v=(1+|x|^2)^{\frac{\alpha}{2}} u^p,
\end{cases} \quad \mbox{in}\;\ \mathbb{R}^N,
\end{equation*}where $N \ge 5$, $p>1$ and $\alpha >0$. We prove Liouville type results for the classical positive (nonnegative) stable solutions in dimension $N<\ell+\dfrac{\alpha (\ell-2)}{2}$ ($N <\ell+\dfrac{\alpha (\ell-2)(p+3)}{4(p+1)}$) and $\ell \ge 5$, $p \in (1,p_*(\ell))$. In particular, for any $p>1$ and $\alpha > 0$, we obtain the nonexistence of classical positive (nonnegative) stable solutions for any $N \le 12+5 \alpha$ ($N\le 12+\dfrac{5\alpha (p+3)}{2(p+1)}$).
\\[0.15cm]
\end{abstract}

\section{Introduction}

We consider the weighted elliptic system
\begin{equation}\label{eq:1.1}
\begin{cases}
-\Delta u=(1+|x|^2)^{\frac{\alpha}{2}} v,\\
-\Delta v=(1+|x|^2)^{\frac{\alpha}{2}} u^p,
\end{cases} \quad \mbox{in}\;\ \mathbb{R}^N,
\end{equation}where $N \ge 5$, $p>1$ and $\alpha>0$. We are interested in the Liouville-type theorems---i.e., the nonexistence of the classical positive and nonnegative stable solutions (\ref{eq:1.1}) in $\mathbb{R}^N$ or the half space $\mathbb{R}^N_+$.\vskip .05in

We recall the case $\alpha=0$, the so-called Lane-Emden equation or system which has been widely studied by many authors.
For the second order Lane-Emden equation, the finite Morse index solutions of the nonlinear problem
\begin{equation}\label{eq:1.2}
\Delta u+|u|^{p-1}u=0\quad \mbox{in}\; \mathbb{R}^N,\; p>1
\end{equation}have been completely classified by Farina (see \cite{Farina}). Farina also proved that nontrivial finite Morse index solutions to (\ref{eq:1.2}) exist if and only if $p \ge p_{JL}$ and $N \ge 11$, or $p=\dfrac{N+2}{N-2}$ and $N \ge 3$. Here $p_{JL}$ is the so-called Joseph-Lundgren exponent (see \cite{Gui}). His proof made a delicate application of the classical Moser's iteration. There exist many excellent papers to utilize Farina's approach to discuss the second order Hardy-H\'{e}non equation. We refer to \cite{Dancer,Wang} and the references therein.\vskip .05in

Unfortunately, Farina's approach may fail to obtain the similarly complete
classification for stable solution and finite Morse index solution of the biharmonic equation
\begin{equation}\label{eq:1.3}
\Delta^2 u=u^p,\quad \mbox{in}\;\ \Omega \subset \mathbb{R}^N.
\end{equation}To solve the complete classification, D\'{a}vila-Dupaigne-Wang-Wei \cite{Davila} have derived from a monotonicity formula for solution of (\ref{eq:1.3}) to reduce the nonexistence of nontrivial entire solutions for the problem (\ref{eq:1.3}), to that of nontrivial homogeneous solutions, and gave a complete classification of stable solutions and those of finite Morse index solutions. Adopting the similar method, Hu \cite{Hu-1} obtained a complete classification of stable solutions and finite Morse index solutions of the fourth order H\'{e}non equation $\Delta^2 u=|x|^{\alpha} |u|^{p-1}u$.\vskip .05in

However, it seems that the monotonicity formula approach in \cite{Davila,Hu-1} does not work well with some weighted elliptic systems or negative exponent. There are several new approaches dealing with those elliptic equation or systems. The first approach is the use of the test function, Souplet's inequality \cite{Souplet} and the idea of Cowan-Esposito-Ghoussoub in \cite{Cowan-Espositio-Ghoussoub}. For example, Fazly proved the following result:\vskip .1in

\noindent {\bf Theorem A}
{\it (\cite[Theorem 2.4]{Fazly}) Suppose that $(u,v) \in C^2(\mathbb{R}^N) \times C^2(\mathbb{R}^N)$ is a nonnegative entire stable solution of (1.1) in dimension
\begin{equation}\label{eq:1.4}
N<8+3\alpha+\dfrac{8+4\alpha}{p-1}.
\end{equation}Then $(u,v)$ has the only trivial solution, where $\alpha \ge 0$ and $p>1$.}
\vskip .1in

The second approach, which was obtained by Cowan-Ghoussoub \cite{Cowan-Ghoussoub} and Dupaigne-Ghergu-Goubet-Warnault \cite{Dupaigne} independently, is firstly to derive the following interesting intermediate second order stability criterion: for the stable positive solution to (\ref{eq:1.3}), it holds
\begin{equation*}
\sqrt{p} \int_{\mathbb{R}^N} u^{\frac{p-1}{2}} \zeta^2 dx \le \int_{\mathbb{R}^N} |\nabla \zeta|^2 dx ,\quad \forall \zeta \in C_0^1(\mathbb{R}^N).
\end{equation*}Then this will be carried out through a bootstrap argument which is reminiscent of the classical Moser iteration method. Recently, combining the first and second approaches, the fourth order elliptic equation with positive or negative exponent have been discussed in \cite{Cowan,Guo-Wei,Hajlaoui}.\vskip .05in

For the general system with $\alpha \neq 0$, the Liouville property is less understood and is more delicate to handle than $\alpha=0$.
Moreover, from Theorem A, we note that if $p>3$ and $\alpha < \dfrac{4(p-3)}{3p+1}$, then the space dimension for the Liouville property of the nonnegative stable solution to (\ref{eq:1.1}) is less than $12$. But the study of radial solutions in \cite{Karageorgis} suggests the following {\bf conjecture}:\\[0.1cm]
\hspace*{18pt}{\it A smooth stable solution to (\ref{eq:1.3}) exists if and only if $p \ge p_{JL_4}$ and $N \ge 13$.}\\[0.1cm]
Consequently, Liouville type result for stable solutions of (\ref{eq:1.1}) should hold true for any $N \le 12$, $p>1$ and $\alpha \ge 0$. That is  what we will prove here.
Inspired by the ideas in \cite{Cowan,Guo-Wei,Hajlaoui}, our purpose in this paper is to prove the following Liouville-type theorems of the weighted elliptic system (\ref{eq:1.1}).

\begin{theorem}\label{eq:t1.1}
Suppose that $(u,v)$ is a classical stable solution of the weighted elliptic system (\ref{eq:1.1}) with $u>0$. If $N \ge 5$, $\alpha>0$ and $p>1$ satisfy the following conditions:
\begin{itemize}
\item [\rm (i).] $
N <\ell+\dfrac{\alpha (\ell-2)}{2}$.
\item [\rm (ii).] $p \in (1,p_*(\ell))$, where
\begin{equation*}p_*(\ell)=
\begin{cases}
+\infty,& 5 \le \ell \le \overline{\ell},\\
\dfrac{\ell+2-\sqrt{\ell^2+4-4\sqrt{\ell^2+H^*_{\ell}}}}{\ell-6-\sqrt{\ell^2+4-4\sqrt{\ell^2+H^*_{\ell}}}}, & \ell > \overline{\ell},
\end{cases}
\end{equation*}and $\overline{\ell} \in (12,13)$ is the root of the quartic equation
\begin{equation*}
8(\ell-2)(\ell-4)=H^*_{\ell}:=\dfrac{\ell^2 (\ell-4)^2}{16}+\dfrac{(\ell-2)^2}{2}-1.
\end{equation*}
\end{itemize}Then the system (\ref{eq:1.1}) has the only trivial solution.
\end{theorem}

\begin{remark}\label{eq:r1.1}
\begin{itemize}
\item [\rm (i).] We note that if $5 \le N\le 12+5 \alpha$, then the weighted elliptic system (\ref{eq:1.1}) do not have classical positive stable solution for any $p>1$ and $\alpha\ge 0$.
\item [\rm (ii).] From (\ref{eq:1.4}) and Theorem \ref{eq:t1.1}, we find that the inequality
\begin{equation*}
8+3\alpha+\dfrac{8+4\alpha}{p-1}<\ell+\dfrac{\alpha(\ell-2)}{2}
\end{equation*}holds true, when $8\le \ell \le 12$, $p=+\infty$ or $\ell \ge 13$, $p=p_*(\ell)$.
\item [\rm (iii).] If we denote $\ell:=2+2\mu$ by Remark \ref{eq:r3.1} (iii), then the weighted elliptic system (\ref{eq:1.1}) do not exist classical positive stable solution in
dimension
\begin{equation*}
N<2+(2+\alpha) \mu,
\end{equation*}where $\alpha \ge 0$, $\mu$ is the largest root of the polynomial
\begin{equation}\label{eq:1.5}
H(p,\mu)=\mu^4-\dfrac{32p(p+1)}{(p-1)^2}\mu^2+\dfrac{32p(p+1)(p+3)}{(p-1)^3}\mu-\dfrac{64p(p+1)^2}{(p-1)^4},
\end{equation}for any $p>1$.
\end{itemize}
\end{remark}\vskip .06in

\begin{theorem}\label{eq:t1.2}
Assume that  $(u,v)$ is a classical nonnegative stable solution of the weighted elliptic system (\ref{eq:1.1}), and
$N\ge 5$, $\alpha>0$ and $p>1$ satisfies one of the following conditions:
\begin{itemize}
\item [\rm (i).] $N <\ell+\dfrac{\alpha (\ell-2)(p+3)}{4(p+1)}$ and $p \in (1,p_*(\ell))$. or
\item [\rm (ii).] For any $p>1$, $N<2+2\mu +\dfrac{\alpha(p+3)}{2(p+1)} \mu$, where $\mu$ is the largest root of the polynomial
$H(p,\mu)$ in (\ref{eq:1.5}).
\end{itemize}Then $(u,v)$ must be the trivial solution.
\end{theorem}

\begin{remark}\label{eq:r1.2}
Clearly, if for any $p>1$ and $\alpha \ge 0$, $5 \le N\le 12+\dfrac{5\alpha (p+3)}{2(p+1)}$, then the system (\ref{eq:1.1}) do no have classical nonnegative stable solution.
\end{remark}

\begin{itemize}
\item {\bf Notation}. Here and in the following, we use $B_r(x)$ to denote the open ball on $\mathbb{R}^N$ central at $x$ with radius $r$. We also write $B_r=B_r(0)$. $C$ denotes generic positive constants independent of $u$, which could be changed from one line to another.
\end{itemize}

The organization of the paper is as follows. In section 2, we prove some decay estimate and point-wise estimate for the stable solution of the system (\ref{eq:1.1}). Then we prove Liouville-type theorem for positive stable solution of (\ref{eq:1.1}), that is Theorem \ref{eq:t1.1} in section 3. Adopting the similar approach, we prove Theorem \ref{eq:t1.2} in section 4.\vskip .1in

\section{Preliminaries}

Let $\Omega$ be a subset of $\mathbb{R}^N$ and $f,g \in C^1\left (\mathbb{R}^{N+2},\Omega \right )$. Following Montenegro \cite{Montenegro}, we consider a general elliptic system
\begin{equation*}
\begin{cases}
-\Delta u=f(u,v,x),\\
-\Delta v=g(u,v,x),
\end{cases}\quad x \in \Omega. \eqno (Q_{f,g})
\end{equation*}

\begin{definition}\label{eq:d1.1}
A solution $(u,v)\in C^2(\Omega)\times C^2(\Omega)$ of $(Q_{f,g})$ is said to be {\it stable}, if the eigenvalue problem
\begin{equation*}
\begin{cases}
-\Delta \phi=f_u(u,v,x)\phi+f_v(u,v,x)\psi+\eta \phi,\\
-\Delta \psi =g_u(u,v,x)\phi+g_v(u,v,x)\psi+\eta \psi,
\end{cases}\eqno (E_{f,g})
\end{equation*} has a first positive eigenvalue $\eta>0$, with corresponding positive smooth eigenvalue pair $(\phi,\psi)$.
A solution $(u,v)$ is called {\it semi-stable}, if the first eigenvalue $\eta$ is nonnegative.
\end{definition}

\begin{lemma}\label{eq:l2.1}
Let $(u,v)$ be a classical stable solution of (\ref{eq:1.1}). Then the following two inequalities hold
\begin{equation*}
p\int_{\mathbb{R}^N} (1+|x|^2)^{\frac{\alpha}{2}} u^{p-1} \zeta^2 dx \le \int_{\mathbb{R}^N} \dfrac{|\Delta \zeta|^2}{(1+|x|^2)^{\frac{\alpha}{2}}}dx, \quad \forall \zeta \in H^2 (\mathbb{R}^N),
\end{equation*}and
\begin{equation}\label{eq:2.1}
\sqrt{p}\int_{\mathbb{R}^N} (1+|x|^2)^{\frac{\alpha}{2}} u^{\frac{p-1}{2}} \zeta^2 dx \le \int_{\mathbb{R}^N} |\nabla \zeta|^2 dx,
\end{equation}for all $\zeta \in H^1(\mathbb{R}^N)$.
\end{lemma}

\begin{proof}
Adopting the proof of Lemma 3 and Lemma 7 in \cite{Cowan}, we get the desired results. Therefore we omit the detail.
\end{proof}

Next, we list some decay estimate and point-wise estimate for stable solution which will be useful in the following proofs.

\begin{lemma}(\cite[Lemma 2.2]{Wei-Ye})\label{eq:l2.2}
For any $\zeta, \eta \in C^4(\mathbb{R}^N)$, the identity holds
\begin{equation*}
\Delta \zeta \Delta (\zeta\eta^2)=[\Delta (\zeta \eta)]^2-4(\nabla \zeta \cdot \nabla \eta)^2-\zeta^2(\Delta \eta)^2+2\zeta \Delta \zeta |\nabla \eta|^2-4\zeta \Delta \eta \nabla \zeta \cdot \nabla \eta.
\end{equation*}
\end{lemma}

\begin{lemma}\label{eq:l2.3}
For any $\zeta \in C^4(\mathbb{R}^N)$ and $\eta \in C^4_0 (\mathbb{R}^N)$, we obtain the two identities
\begin{align}\label{eq:2.2}
\int_{\mathbb{R}^N} \dfrac{\Delta \zeta \Delta (\zeta \eta^2)}{(1+|x|^2)^{\frac{\alpha}{2}}} dx =& \int_{\mathbb{R}^N} \dfrac{\left [\Delta (\zeta \eta) \right ]^2}{(1+|x|^2)^{\frac{\alpha}{2}}}dx +\int_{\mathbb{R}^N}
\dfrac{\left [-4(\nabla \zeta \cdot \nabla \eta)^2+2 \zeta \Delta \zeta | \nabla \eta|^2 \right ]}{(1+|x|^2)^{\frac{\alpha}{2}}}dx \nonumber \\[0.1cm]
 & +\int_{\mathbb{R}^N} \dfrac{\zeta^2}{(1+|x|^2)^{\frac{\alpha}{2}}} \left [2\nabla (\Delta \eta)\cdot \nabla \eta+(\Delta \eta)^2 \right ]dx \nonumber \\[0.1cm]
 & -2\alpha \int_{\mathbb{R}^N} \dfrac{\zeta^2}{(1+|x|^2)^{\frac{\alpha}{2}+1}}\Delta \eta (\nabla \eta \cdot x) dx,
\end{align}and
\begin{align}\label{eq:2.3}
2 \int_{\mathbb{R}^N} \dfrac{|\nabla \zeta|^2|\nabla \eta|^2}{(1+|x|^2)^{\frac{\alpha}{2}}}dx= & 2\int_{\mathbb{R}^N} \dfrac{\zeta (-\Delta \zeta)|\nabla \eta|^2 }{(1+|x|^2)^{\frac{\alpha}{2}}}dx+\int_{\mathbb{R}^N} \dfrac{\zeta^2 \Delta (|\nabla \eta|^2) }{(1+|x|^2)^{\frac{\alpha}{2}}}dx \nonumber \\[0.1cm]
& -\alpha \int_{\mathbb{R}^N} \dfrac{\zeta^2}{(1+|x|^2)^{\frac{\alpha}{2}+1}}\left [2(\nabla (|\nabla \eta|^2) \cdot x)+N|\nabla \eta|^2\right ] dx \nonumber \\[0.1cm]
&+\alpha(\alpha+2)\int_{\mathbb{R}^N} \dfrac{\zeta^2 |x|^2 |\nabla \eta|^2}{(1+|x|^2)^{\frac{\alpha}{2}+2}}dx.
\end{align}
\end{lemma}

\begin{proof}
Integrating by parts, we get
\begin{align*}
-4 \int_{\mathbb{R}^N} & \dfrac{\zeta \Delta \eta \nabla \zeta \cdot \nabla \eta }{(1+|x|^2)^{\frac{\alpha}{2}}}dx
= 2\int_{\mathbb{R}^N}\zeta^2\cdot div \left ( \dfrac{\Delta \eta \nabla \eta }{(1+|x|^2)^{\frac{\alpha}{2}}} \right )dx\\[0.1cm]
= & 2 \int_{\mathbb{R}^N} \dfrac{\zeta^2}{(1+|x|^2)^{\frac{\alpha}{2}}}\left [\nabla (\Delta \eta)\cdot \nabla \eta +(\Delta \eta)^2 \right ] dx-2\alpha \int_{\mathbb{R}^N} \dfrac{\zeta^2 \Delta \eta (\nabla \eta \cdot x)}{(1+|x|^2)^{\frac{\alpha}{2}+1}}dx.
\end{align*}Combining with Lemma \ref{eq:l2.2}, it implies that the identity (\ref{eq:2.2}) holds true.

A simple computation leads to
\begin{align*}
\int_{\mathbb{R}^N} \dfrac{\Delta (\zeta^2)|\nabla \eta|^2 }{(1+|x|^2)^{\frac{\alpha}{2}}}dx = & \int_{\mathbb{R}^N} \zeta^2 \Delta \left ( \dfrac{ |\nabla \eta|^2 }{(1+|x|^2)^{\frac{\alpha}{2}}} \right )dx\\[0.1cm]
=& \int_{\mathbb{R}^N} \dfrac{\zeta^2 \Delta (|\nabla \eta|^2) }{(1+|x|^2)^{\frac{\alpha}{2}}}dx-\alpha \int_{\mathbb{R}^N} \dfrac{\zeta^2 \left [2(\nabla (|\nabla \eta|^2) \cdot x)+N|\nabla \eta|^2\right ]}{(1+|x|^2)^{\frac{\alpha}{2}+1}} dx \\[0.1cm]
&+\alpha(\alpha+2)\int_{\mathbb{R}^N} \dfrac{\zeta^2 |x|^2 |\nabla \eta|^2}{(1+|x|^2)^{\frac{\alpha}{2}+2}}dx.
\end{align*}Again it is easy to verify that
\begin{equation*}
\dfrac{1}{2}\Delta (\zeta^2)=\zeta \Delta \zeta+|\nabla \zeta|^2.
\end{equation*}Combining the above two identities, we get the identity (\ref{eq:2.3}).
\end{proof}

\begin{lemma}\label{eq:l2.4}
Let $N \ge 5$, $p>1$, $\alpha \ge 0$ and $(u,v)$ be a classical stable solution of (\ref{eq:1.1}) with $u \ge 0$. Then we have
\begin{equation*}
\int_{B_R} (1+|x|^2)^{\frac{\alpha}{2}} \left [v^2+u^{p+1} \right ] dx \le CR^{N-4-\alpha-\frac{8+4\alpha}{p-1}},
\end{equation*}for all $R>0$.
\end{lemma}

\begin{proof}
Since $(u,v)$ is a classical stable solution of (\ref{eq:1.1}), we find that for any $\zeta \in C_0^4(\mathbb{R}^N)$,
\begin{equation}\label{eq:2.4}
\int_{\mathbb{R}^N} (1+|x|^2)^{\frac{\alpha}{2}} u^p \zeta dx=\int_{\mathbb{R}^N} \dfrac{\Delta u}{(1+|x|^2)^{\frac{\alpha}{2}}} \Delta \zeta dx,
\end{equation}and
\begin{equation}\label{eq:2.5}
p\int_{\mathbb{R}^N} (1+|x|^2)^{\frac{\alpha}{2}}u^{p-1}\zeta^2 dx \le \int_{\mathbb{R}^N} \dfrac{|\Delta \zeta|^2}{(1+|x|^2)^{\frac{\alpha}{2}}}dx.
\end{equation}Substituting $\zeta =u \psi^2$ into (\ref{eq:2.4}), we obtain
\begin{equation}\label{eq:2.6}
\int_{\mathbb{R}^N} (1+|x|^2)^{\frac{\alpha}{2}} u^{p+1}\psi^2 dx=\int_{\mathbb{R}^N} \dfrac{\Delta u \Delta (u\psi^2)}{(1+|x|^2)^{\frac{\alpha}{2}}} dx.
\end{equation}Substitute $\zeta=u\psi$ into (\ref{eq:2.5}) to get
\begin{equation}\label{eq:2.7}
p\int_{\mathbb{R}^N} (1+|x|^2)^{\frac{\alpha}{2}}u^{p+1}\psi^2 dx \le \int_{\mathbb{R}^N} \dfrac{[\Delta (u\psi)]^2}{(1+|x|^2)^{\frac{\alpha}{2}}}dx.
\end{equation}Here and in the following, we choose the cut-off function $\psi \in C^4_0(\mathbb{R}^N)$ with $0 \le \psi \le 1$,
\begin{equation*}
\psi (x)=
\begin{cases}
1, & \mbox{if}\;\ |x|<R, \\\
0, & \mbox{if}\;\ |x|>2R,
\end{cases}
\end{equation*}and $|\nabla^i \psi |\le \dfrac{C}{R^i}$, for $i=1,2,3$. Now, combining (\ref{eq:2.6}) and (\ref{eq:2.7}) with (\ref{eq:2.2}) and (\ref{eq:2.3}), we have
\begin{align*}
(p-1)\int_{\mathbb{R}^N} & (1+|x|^2)^{\frac{\alpha}{2}} u^{p+1}\psi^2 dx \le \int_{\mathbb{R}^N} \dfrac{4(\nabla u\cdot \nabla \psi)^2-2u\Delta u|\nabla \psi|^2}{(1+|x|^2)^{\frac{\alpha}{2}}}dx\\[0.12cm]
& -\int_{\mathbb{R}^N} \dfrac{u^2 \left [2\nabla (\Delta \psi)\cdot \nabla \psi+|\Delta \psi|^2\right ] }{(1+|x|^2)^{\frac{\alpha}{2}}}dx +2\alpha \int_{\mathbb{R}^N} \dfrac{u^2 \Delta \psi (\nabla \psi \cdot x)}{(1+|x|^2)^{\frac{\alpha}{2}+1}} dx\\[0.12cm]
\le & C\int_{\mathbb{R}^N} |uv| \cdot |\nabla \psi|^2 dx\\[0.1cm]
& +C\int_{\mathbb{R}^N} \dfrac{u^2}{(1+|x|^2)^{\frac{\alpha}{2}}}\left [|\Delta (|\nabla \psi|^2)|+|\nabla (\Delta \psi)\cdot \nabla \psi|+|\Delta \psi|^2\right ] dx\\[0.1cm]
&+C \int_{\mathbb{R}^N} \dfrac{u^2}{(1+|x|^2)^{\frac{\alpha}{2}+1}}\left [|\nabla (|\nabla \psi|^2)\cdot x|+|\nabla \psi|^2+|\Delta \psi (\nabla \psi \cdot x)| \right ]dx.
\end{align*}Again since $\Delta (u\psi)=-(1+|x|^2)^{\frac{\alpha}{2}}v\psi+2\nabla u\cdot \nabla \psi+u\Delta \psi$, we get
\begin{equation*}
\int_{\mathbb{R}^N}(1+|x|^2)^{\frac{\alpha}{2}}v^2\psi^2 dx \le C\int_{\mathbb{R}^N} \dfrac{1}{(1+|x|^2)^{\frac{\alpha}{2}}}\left [(\Delta (u\psi))^2+|\nabla u\cdot \nabla \psi|^2+u^2|\Delta \psi|^2 \right ] dx.
\end{equation*}Then, combining the above inequality with (\ref{eq:2.6}), it implies that
\begin{align*}
\int_{\mathbb{R}^N} & (1+|x|^2)^{\frac{\alpha}{2}}\left [v^2+u^{p+1}\right ]\psi^2 dx \le C \int_{\mathbb{R}^N} |uv| |\nabla \psi|^2dx \\
& +C \int_{\mathbb{R}^N} \dfrac{u^2}{(1+|x|^2)^{\frac{\alpha}{2}}}\left [|\Delta (|\nabla \psi|^2)|+|\nabla (\Delta \psi)\cdot \nabla \psi|+|\Delta \psi|^2\right ]dx\\[0.1cm]
& +C \int_{\mathbb{R}^N} \dfrac{u^2}{(1+|x|^2)^{\frac{\alpha}{2}+1}}\left [|\nabla (|\nabla \psi|^2)\cdot x|+|\nabla \psi|^2+|\Delta \psi (\nabla \psi \cdot x)| \right ]dx.
\end{align*}

Next, the function $\psi$ in the above inequality are replaced by $\psi^m$, where $m$ is a larger integer, then
\begin{align}\label{eq:2.8}
\int_{\mathbb{R}^N} & (1+|x|^2)^{\frac{\alpha}{2}} \left [v^2+u^{p+1} \right ]\psi^{2m} dx \le C \int_{\mathbb{R}^N} |uv| \psi^{2(m-1)}|\nabla \psi|^2dx \nonumber \\[0.1cm]
& +C\int_{\mathbb{R}^N} \dfrac{u^2}{(1+|x|^2)^{\frac{\alpha}{2}}}\left [|\Delta (|\nabla \psi^m|^2)|+|\nabla (\Delta \psi^m)\cdot \nabla \psi^m|+|\Delta \psi^m|^2 \right ]dx \nonumber \\[0.1cm]
& +C \int_{\mathbb{R}^N} \dfrac{u^2}{(1+|x|^2)^{\frac{\alpha}{2}+1}}\left [|\nabla (|\nabla \psi^m|^2)\cdot x|+|\nabla \psi^m|^2+|\Delta \psi^m (\nabla \psi^m \cdot x)| \right ]dx.
\end{align}A simple application of Young's inequality leads to
\begin{align*}
\int_{\mathbb{R}^N} |uv|\psi^{2(m-1)}|\nabla \psi|^2 dx \le \dfrac{1}{2C} \int_{\mathbb{R}^N}(1+|x|^2)^{\frac{\alpha}{2}}v^2 \psi^{2m}dx\\[0.1cm]
+C\int_{\mathbb{R}^N} \dfrac{u^2}{(1+|x|^2)^{\frac{\alpha}{2}}} \psi^{2(m-2)} |\nabla \psi|^4 dx,
\end{align*}and putting into (\ref{eq:2.8}) yields
\begin{align*}
\int_{\mathbb{R}^N} & (1+|x|^2)^{\frac{\alpha}{2}}\left [v^2+u^{p+1}\right ]\psi^{2m} dx \le C \int_{\mathbb{R}^N} \dfrac{u^2}{(1+|x|^2)^{\frac{\alpha}{2}}} \psi^{2(m-2)}|\nabla \psi|^4 dx \\[0.1cm]
& +C\int_{\mathbb{R}^N} \dfrac{u^2}{(1+|x|^2)^{\frac{\alpha}{2}}}\left [|\Delta (|\nabla \psi^m|^2)|+|\nabla (\Delta \psi^m)\cdot \nabla \psi^m|+|\Delta \psi^m|^2\right ]dx\\[0.1cm]
& +C \int_{\mathbb{R}^N} \dfrac{u^2}{(1+|x|^2)^{\frac{\alpha}{2}+1}}\left [|\nabla (|\nabla \psi^m|^2)\cdot x|+|\nabla \psi^m|^2+|\Delta \psi^m (\nabla \psi^m \cdot x)|\right ]dx\\[0.1cm]
=&C\int_{\mathbb{R}^N} \dfrac{u^2}{(1+|x|^2)^{\frac{\alpha}{2}}}\psi^{2(m-2)} \mathfrak{B}(\psi^m) dx
+C\int_{\mathbb{R}^N} \dfrac{u^2}{(1+|x|^2)^{\frac{\alpha}{2}+1}}\psi^{2(m-2)} \mathfrak{F}(\psi^m) dx,
\end{align*}where
\begin{align*}
\mathfrak{B}(\psi^m) & =|\nabla \psi|^4 +\psi^{2(2-m)}\left [|\Delta (|\nabla \psi^m|^2)|+|\nabla (\Delta \psi^m)\cdot \nabla \psi^m|+|\Delta \psi^m|^2\right ],\\[0.1cm]
\mathfrak{F}(\psi^m) & =\psi^{2(2-m)} \left [|\nabla (|\nabla \psi^m|^2)\cdot x|+|\nabla \psi^m|^2+|\Delta \psi^m (\nabla \psi^m \cdot x)|\right ].
\end{align*}Choosing $m$ larger enough such that $(m-2)(p+1) \ge 2m$, we utilize H\"{o}lder's inequality to the both terms in the right side of the above inequality and find
\begin{align*}
\int_{\mathbb{R}^N} & \dfrac{u^2}{(1+|x|^2)^{\frac{\alpha}{2}}}\psi^{2(m-2)} \mathfrak{B}(\psi^m) dx \\[0.1cm]
= & \int_{\mathbb{R}^N} (1+|x|^2)^{\frac{\alpha}{p+1}}u^2\psi^{2(m-2)} (1+|x|^2)^{-\frac{\alpha}{p+1}-\frac{\alpha}{2}} \mathfrak{B}(\psi^m)dx \\[0.1cm]
\le & \left (\int_{\mathbb{R}^N} (1+|x|^2)^{\frac{\alpha}{2}}u^{p+1}\psi^{2m} dx \right )^{\frac{2}{p+1}}\\[0.1cm]
& \times \left ( \int_{\mathbb{R}^N} (1+|x|^2)^{-\frac{2\alpha+\alpha(p+1)}{2(p-1)}} \mathfrak{B}(\psi^m)^{\frac{p+1}{p-1}} dx \right )^{\frac{p-1}{p+1}}
\end{align*}and
\begin{align*}
\int_{\mathbb{R}^N} & \dfrac{u^2}{(1+|x|^2)^{\frac{\alpha}{2}+1}}\psi^{2(m-2)} \mathfrak{F}(\psi^m) dx\\[0.1cm]
\le & \left (\int_{\mathbb{R}^N} (1+|x|^2)^{\frac{\alpha}{2}}u^{p+1}\psi^{2m} dx \right )^{\frac{2}{p+1}} \\[0.1cm]
&\times \left ( \int_{\mathbb{R}^N} (1+|x|^2)^{-\frac{2\alpha+(\alpha+2)(p+1)}{2(p-1)}} \mathfrak{F}(\psi^m)^{\frac{p+1}{p-1}} dx \right )^{\frac{p-1}{p+1}}.
\end{align*}Therefore, we get
\begin{align*}
\int_{\mathbb{R}^N} & (1+|x|^2)^{\frac{\alpha}{2}}\left [ v^2+u^{p+1} \right ] \psi^{2m} dx \\[0.1cm]
\le & C\int_{\mathbb{R}^N} (1+|x|^2)^{-\frac{2\alpha+\alpha(p+1)}{2(p-1)}} \mathfrak{B}(\psi^m)^{\frac{p+1}{p-1}} dx\\[0.1cm]
&+C \int_{\mathbb{R}^N} (1+|x|^2)^{-\frac{2\alpha+(\alpha+2)(p+1)}{2(p-1)}} \mathfrak{F}(\psi^m)^{\frac{p+1}{p-1}} dx\\[0.1cm]
\le & C R^{N-4-\alpha-\frac{8+4\alpha}{p-1}},
\end{align*}for all $R>0$.
\end{proof}

Let $N\ge 5$, $p>1$ and $\alpha >0$. We consider a more general elliptic system
\begin{equation}\label{eq:2.9}
\begin{cases}
\begin{split}
&-\Delta u=(1+|x|^2)^{\frac{\alpha}{2}}v,\\
&-\Delta v=(1+|x|^2)^{\frac{\alpha}{2}}u^p,
\end{split}\quad & \mbox{in}\; \Sigma, \\
\ u=\Delta u=0, \quad & \mbox{on}\; \partial \Sigma,
\end{cases}
\end{equation}where $\Sigma=\mathbb{R}^N$ or the half space $\Sigma=\mathbb{R}^N_+$ or the exterior domain $\Sigma=\mathbb{R}^N \backslash \overline{\Omega}$, $\mathbb{R}^N_+\backslash \overline{\Omega}$, and $\Omega$ is a bounded smooth domain of $\mathbb{R}^N$. A solution $(u,v)$ of (\ref{eq:2.9}) is said to be stable if for any $\zeta \in H^2(\Sigma) \cap H_0^1(\Sigma)$, we have
\begin{equation*}
p\int_{\Sigma} (1+|x|^2)^{\frac{\alpha}{2}} u^{p-1} \zeta^2 dx \le \int_{\Sigma} \dfrac{|\Delta \zeta|^2}{(1+|x|^2)^{\frac{\alpha}{2}}}dx,
\end{equation*}or if for any $\zeta \in H^1(\Sigma)$
\begin{equation*}
\sqrt{p}\int_{\Sigma} (1+|x|^2)^{\frac{\alpha}{2}} u^{\frac{p-1}{2}} \zeta^2 dx \le \int_{\Sigma} |\nabla \zeta|^2 dx.
\end{equation*}

Motivated by the proof in \cite{Hajlaoui,Phan,Souplet}, we obtain the crucial ingredient in the proof of Theorem \ref{eq:t1.1} and Theorem \ref{eq:t1.2}.

\begin{lemma}\label{eq:l2.5}
Assume that $(u,v)$ is a classical stable solution of (\ref{eq:2.9}). If $u\ge 0$, then the inequality holds
\begin{equation*}
v \ge \sqrt{\dfrac{2}{p+1}}u^{\frac{p+1}{2}},\quad \mbox{in}\;\ \Sigma.
\end{equation*}
\end{lemma}

\begin{proof}
Set $\delta=\sqrt{\dfrac{2}{p+1}}$, $\gamma=\delta u^{\frac{p+1}{2}} -v$. A direct calculation yields
\begin{align*}
\Delta \gamma & =\frac{p-1}{2\delta}\cdot u^{\frac{p-3}{2}}|\nabla u|^2-\delta^{-1}(1+|x|^2)^{\frac{\alpha}{2}}u^{\frac{p-1}{2}}v+(1+|x|^2)^{\frac{\alpha}{2}}u^p \\
& \ge  \delta^{-1}(1+|x|^2)^{\frac{\alpha}{2}}u^{\frac{p-1}{2}}\gamma.
\end{align*}Denote $\gamma_+:=\max \{\gamma,0\}$. Since $\gamma_+\Delta \gamma \ge 0$ in $\Sigma$ and $\gamma=0$ on $\partial \Sigma$,
 it implies that for any $R>0$
\begin{align}\label{eq:2.10}
\int_{\Sigma \cap B_R} |\nabla \gamma_+|^2dx & =-\int_{\Sigma \cap B_R}\gamma_+\Delta \gamma dx+\int_{\partial (\Sigma \cap B_R)}\gamma_+\dfrac{\partial \gamma}{\partial \nu} d\sigma \nonumber \\[0.1cm]
& \le \int_{\Sigma \cap \partial B_R} \gamma_+\dfrac{\partial \gamma}{\partial \nu}d \sigma.
\end{align}For $r>0$, define $\xi (r):=\displaystyle \int_{\mathbb{S}^{N-1} \cap (r^{-1} \Sigma)}\gamma_+^2(r\sigma)d \sigma$, where $\mathbb{S}^{N-1}$ denotes by the unit sphere in $\mathbb{R}^N$.  We easily deduce that there exists $R_0>0$ such that
\begin{equation*}
\int_{\Sigma \cap \partial B_r} \gamma_+\dfrac{\partial \gamma}{\partial \nu} d \sigma =\dfrac{r^{N-1}}{2}\xi '(r),\quad \forall r \ge R_0.
\end{equation*}

On the other hand, for any $R \ge R_0$, we conclude that
\begin{align*}
\int_{R_0}^R r^{N-1} \xi (r)dr = & \int_{R_0}^R \int_{\mathbb{S}^{N-1}\cap (r^{-1}\Sigma)} \gamma_+^2(r\sigma) r^{N-1} dr d\sigma \\[0.1cm]
\le & C\int_{B_R\cap \Sigma}\gamma_+^2 dx \le C\int_{B_R\cap \Sigma}(1+|x|^2)^{\frac{\alpha}{2}}(v^2+u^{p+1})dx\\[0.1cm]
\le & CR^{N-4-\alpha-\frac{8+4\alpha}{p-1}}=o(R^N).
\end{align*}Therefore, $\xi (R_i) \to 0$ for some sequence $R_i \to \infty$. Thus, there exists $\tilde{R}_i \to +\infty$ such that $\xi'(\tilde{R}_i)\le 0$. Letting $i \to \infty$ in (\ref{eq:2.10}) with $R=\tilde{R}_i$, we find
\begin{equation*}
\int_{\Sigma}|\nabla \gamma_+|^2 dx =0.
\end{equation*}Again sine $\gamma=0$ on $\partial \Sigma$, we get that $\gamma_+\equiv 0$ in $\Sigma$.
\end{proof}

Throughout the paper, we let $R_k=2^k R$ with $R>0$ and integers $k \ge 1$.

\begin{lemma}(\cite[Lemma 5]{Cowan})\label{eq:l2.6}
For any integer $k\ge 1$ and $1 \le \beta <\dfrac{N}{N-2}$, there is some $C=C(k,\beta)<+\infty$ such that for any smooth $w \ge 0$, the inequality holds
\begin{equation*}
\left (\int_{B_{R_k}} w^{\beta}dx \right )^{\frac{1}{\beta}}\le CR^{2+N(\frac{1}{\beta}-1)}\int_{B_{R_{k+1}}}|\Delta w|dx +CR^{N(\frac{1}{\beta}-1)}
\int_{B_{R_{k+1}}}wdx.
\end{equation*}
\end{lemma}\vskip .3in

\section{Proof of Theorem \ref{eq:t1.1}}

The following two lemmas play an important role in dealing
with Theorem \ref{eq:t1.1}.

\begin{lemma}\label{eq:l3.1}
Let $N \ge 5$, $p>1$ and $\alpha >0$.
Assume that $(u,v)$ is a classical nonnegative stable solution of (\ref{eq:1.1}). Then we obtain that, for any $s>2$ and $R>0$
\begin{equation*}
\int_{B_{R_k}} (1+|x|^2)^{\frac{\alpha}{2}} u^p v^{s-1} dx \le \dfrac{C}{R^2} \int_{B_{R_{k+1}}} v^s dx,
\end{equation*}provided that
\begin{equation}\label{eq:3.1}
L(p,s)=s^4-\dfrac{32p}{p+1}s^2+\dfrac{32p(p+3)}{(p+1)^2}s-\dfrac{64p}{(p+1)^2} <0.
\end{equation}
\end{lemma}

\begin{proof}
Testing (\ref{eq:2.1}) on $\zeta=u^{\frac{q+1}{2}}\phi$ with $\phi \in C^2_0(\mathbb{R}^N)$ and $ q \ge 1$,
we get
\begin{align}\label{eq:3.2}
\sqrt{p}\int_{\mathbb{R}^N} & (1+|x|^2)^{\frac{\alpha}{2}} u^{\frac{p-1}{2}}u^{q+1}\phi^2 dx \le \int_{\mathbb{R}^N} u^{q+1}|\nabla \phi|^2 dx \nonumber \\[0.1cm]
& +\int_{\mathbb{R}^N} |\nabla u^{\frac{q+1}{2}}|^2 \phi^2 dx+(q+1)\int_{\mathbb{R}^N} u^q \phi \nabla u\cdot \nabla \phi dx.
\end{align}Integrating by parts, we have
\begin{align*}
(q+1)\int_{\mathbb{R}^N} u^q \phi & \nabla u\cdot \nabla \phi dx= \dfrac{1}{2} \int_{\mathbb{R}^N} \nabla (u^{q+1})\cdot \nabla (\phi^2) dx \\[0.1cm]
= & -\dfrac{1}{2}\int_{\mathbb{R}^N} u^{q+1}\Delta (\phi^2) dx,
\end{align*}and
\begin{align*}
\int_{\mathbb{R}^N} |\nabla u^{\frac{q+1}{2}}|^2 & \phi^2 dx = \dfrac{(q+1)^2}{4q}\int_{\mathbb{R}^N} \phi^2 \nabla (u^q)\cdot \nabla u dx \\[0.1cm]
=& -\dfrac{(q+1)^2}{4q}\int_{\mathbb{R}^N} u^q\phi^2 \Delta u dx -\dfrac{q+1}{4q}\int_{\mathbb{R}^N} \nabla (u^{q+1})\nabla (\phi^2) dx\\[0.1cm]
=& \dfrac{(q+1)^2}{4q}\int_{\mathbb{R}^N}(1+|x|^2)^{\frac{\alpha}{2}} u^q v \phi^2 dx+\dfrac{q+1}{4q}\int_{\mathbb{R}^N} u^{q+1} \Delta (\phi^2) dx.
\end{align*}Combining the above two identities with (\ref{eq:3.2}) leads to
\begin{align*}
a_1\int_{\mathbb{R}^N} (1+|x|^2)^{\frac{\alpha}{2}} & u^{\frac{p-1}{2}}u^{q+1}\phi^2 dx\le \int_{\mathbb{R}^N} (1+|x|^2)^{\frac{\alpha}{2}} u^qv\phi^2 dx\\[0.1cm]
& +C\int_{\mathbb{R}^N} u^{q+1}\left [|\Delta (\phi^2)| +|\nabla \phi|^2\right ]dx,
\end{align*}where $a_1=\dfrac{4q\sqrt{p}}{(q+1)^2}$.
We choose $\phi(x)=\omega \left (\dfrac{x}{R_k} \right )$ where $\omega \in C_0^2 (B_2)$ with $\omega \equiv 1$ in $B_1$, then we find
\begin{equation}\label{eq:3.3}
\int_{\mathbb{R}^N} (1+|x|^2)^{\frac{\alpha}{2}} u^{\frac{p-1}{2}}u^{q+1}\phi^2 dx\le \dfrac{1}{a_1}\int_{\mathbb{R}^N} (1+|x|^2)^{\frac{\alpha}{2}} u^qv\phi^2 dx +\dfrac{C}{R^2}\int_{B_{R_{k+1}}}
 u^{q+1}dx.
\end{equation}Similarly, testing (\ref{eq:2.1}) on $\zeta=v^{\frac{r+1}{2}}\phi$ with $r \ge 1$, there holds
\begin{align*}
\sqrt{p}\int_{\mathbb{R}^N} & (1+|x|^2)^{\frac{\alpha}{2}} u^{\frac{p-1}{2}}v^{r+1}\phi^2 dx \le \int_{\mathbb{R}^N} v^{r+1}|\nabla \phi|^2 dx\\[0.1cm]
& +\int_{\mathbb{R}^N} |\nabla v^{\frac{r+1}{2}}|^2 \phi^2 dx+(r+1)\int_{\mathbb{R}^N} v^r\phi \nabla v\cdot \nabla \phi dx.
\end{align*}Adopting the same computation as above (noting that the equation $-\Delta v=(1+|x|^2)^{\frac{\alpha}{2}}u^p$), we obtain
\begin{equation}\label{eq:3.4}
\int_{\mathbb{R}^N} (1+|x|^2)^{\frac{\alpha}{2}} u^{\frac{p-1}{2}}v^{r+1}\phi^2 dx\le \dfrac{1}{a_2}\int_{\mathbb{R}^N} (1+|x|^2)^{\frac{\alpha}{2}} u^pv^r\phi^2 dx
+\dfrac{C}{R^2}\int_{B_{R_{k+1}}}
 v^{r+1}dx,
\end{equation}where $a_2=\dfrac{4r\sqrt{p}}{(r+1)^2}$. \vskip .1in

Rewriting (\ref{eq:3.3}) and (\ref{eq:3.4}) yields
\begin{align}\label{eq:3.5}
I_1+a_2^{r+1}I_2:= & \int_{\mathbb{R}^N} (1+|x|^2)^{\frac{\alpha}{2}} u^{\frac{p-1}{2}}u^{q+1}\phi^2 dx +a_2^{r+1}\int_{\mathbb{R}^N} (1+|x|^2)^{\frac{\alpha}{2}} u^{\frac{p-1}{2}}v^{r+1} \phi^2 dx \nonumber \\[0.1cm]
\le & \dfrac{1}{a_1}\int_{\mathbb{R}^N}(1+|x|^2)^{\frac{\alpha}{2}} u^q v\phi^2 dx +a_2^r\int_{\mathbb{R}^N}(1+|x|^2)^{\frac{\alpha}{2}}u^p v^r \phi^2 dx \nonumber \\[0.1cm]
& +\dfrac{C}{R^2}\int_{B_{R_{k+1}}} (u^{q+1}+v^{r+1}) dx.
\end{align}Fix now
\begin{equation*}
2q=(p+1)r+p-1 \Longleftrightarrow q+1=\dfrac{(p+1)(r+1)}{2}.
\end{equation*}A direct application of Young's inequality leads to
\begin{align*}
\dfrac{1}{a_1}\int_{\mathbb{R}^N} (1+|x|^2 & )^{\frac{\alpha}{2}} u^q v\phi^2 dx = \int_{\mathbb{R}^N} (1+|x|^2)^{\frac{\alpha}{2}} u^{\frac{p-1}{2}}\phi^2 u^{\frac{(q+1)r}{r+1}}\dfrac{v}{a_1} dx\\[0.08cm]
\le & \dfrac{r}{r+1} \int_{\mathbb{R}^N} (1+|x|^2)^{\frac{\alpha}{2}}u^{\frac{p-1}{2}}u^{q+1}\phi^2 dx \\[0.08cm]
& +\dfrac{1}{a_1^{r+1}(r+1)}\int_{\mathbb{R}^N} (1+|x|^2)^{\frac{\alpha}{2}} u^{\frac{p-1}{2}}v^{r+1}\phi^2 dx\\[0.08cm]
= & \dfrac{r}{r+1}I_1+\dfrac{1}{a_1^{r+1}(r+1)} I_2.
\end{align*}Arguing as above, we get
\begin{align*}
a_2^r\int_{\mathbb{R}^N} (1+|x|^2)^{\frac{\alpha}{2}} & u^pv^r \phi^2 dx
\le \dfrac{1}{r+1}\int_{\mathbb{R}^N}(1+|x|^2)^{\frac{\alpha}{2}}u^{\frac{p-1}{2}}u^{q+1}\phi^2 dx \\[0.08cm]
& +\dfrac{a_2^{r+1}r}{r+1}\int_{\mathbb{R}^N} (1+|x|^2)^{\frac{\alpha}{2}}u^{\frac{p-1}{2}}v^{r+1}\phi^2 dx\\[0.08cm]
= & \dfrac{1}{r+1}I_1+\dfrac{a_2^{r+1}r}{r+1}I_2.
\end{align*}Combining the above two inequalities with (\ref{eq:3.5}), we find
\begin{align*}
I_1+a_2^{r+1}I_2\le & \dfrac{r}{r+1}I_1+\dfrac{1}{a_1^{r+1}(r+1)}I_2+\dfrac{1}{r+1}I_1+\dfrac{a_2^{r+1}r}{r+1}I_2\\[0.08cm]
& +\dfrac{C}{R^2}\int_{B_{R_{k+1}}} (u^{q+1}+v^{r+1})dx.
\end{align*}Thus, it implies that
\begin{equation*}
\dfrac{(a_1a_2)^{r+1}-1}{r+1}I_2 \le \dfrac{Ca_1^{r+1}}{R^2}\int_{B_{R_{k+1}}} (u^{q+1}+v^{r+1})dx.
\end{equation*}If $a_1a_2>1$, then we conclude from the choice of $\phi$ that
\begin{equation*}
\int_{B_{R_k}} (1+|x|^2)^{\frac{\alpha}{2}}u^{\frac{p-1}{2}}v^{r+1} dx \le I_2 \le \dfrac{C}{R^2} \int_{B_{R_{k+1}}} (u^{q+1}+v^{r+1})dx.
\end{equation*}Since $q+1=\dfrac{(p+1)(r+1)}{2}$ and $v \ge \sqrt{\dfrac{2}{p+1}}u^{\frac{p+1}{2}}$, we have
\begin{equation*}
\int_{B_{R_k}} (1+|x|^2)^{\frac{\alpha}{2}}u^{\frac{p-1}{2}}v^{r+1}dx \le \dfrac{C}{R^2}\int_{B_{R_{k+1}}} v^{r+1}dx.
\end{equation*}Denote $s:=r+1$. Then, combining the above inequality with Lemma \ref{eq:l2.5}, it implies that
\begin{align*}
\int_{B_{R_k}} (1+|x|^2)^{\frac{\alpha}{2}}u^pv^{s-1} dx =\int_{B_{R_k}} (1+|x|^2)^{\frac{\alpha}{2}}u^{\frac{p-1}{2}}u^{\frac{p+1}{2}}v^{s-1}dx\\[0.08cm]
\le \int_{B_{R_k}} (1+|x|^2)^{\frac{\alpha}{2}} u^{\frac{p-1}{2}}v^s dx \le \dfrac{C}{R^2}\int_{B_{R_{k+1}}} v^sdx.
\end{align*}

On the other hand, since $a_1=\dfrac{4q\sqrt{p}}{(q+1)^2}$, $a_2=\dfrac{4r\sqrt{p}}{(r+1)^2}$, and $2q=(p+1)r+p-1$, $q+1=\dfrac{p+1}{2}(r+1)$, $s=r+1$, we obtain that
\begin{equation*}
a_1a_2>1
\end{equation*}is equivalent to
\begin{equation*}
8p[(p+1)(s-1)+p-1](s-1) >\dfrac{(p+1)^2}{4}s^4.
\end{equation*}A direct calculation shows
\begin{equation*}
L(p,s):=s^4-\dfrac{32p}{p+1}s^2+\dfrac{32p(p+3)}{(p+1)^2}s-\dfrac{64p}{(p+1)^2}<0.
\end{equation*}This completes the proof.
\end{proof}\vskip .1in

\begin{remark}\label{eq:r3.1}
For any $p>1$, the following statements are equivalent:
\begin{itemize}
\item [\rm (i).] $L(p,s)=s^4-\dfrac{32p}{p+1}s^2+\dfrac{32p(p+3)}{(p+1)^2}s-\dfrac{64p}{(p+1)^2}<0$.
\item [\rm (ii).] $H(p,\mu) =\mu^4-\dfrac{32p(p+1)}{(p-1)^2}\mu^2+\dfrac{32p(p+1)(p+3)}{(p-1)^3}\mu-\dfrac{64p(p+1)^2}{(p-1)^4}<0$.
\item [\rm (iii).] Let $\ell=2+2\mu$, then we have $p\in (1,p_*(\ell))$.
\end{itemize}
\end{remark}

\begin{proof}
(i) $\Longleftrightarrow$ (ii). Set $\mu: =\dfrac{p+1}{p-1}s$. A simple calculation yields
\begin{align}\label{eq:3.6}
\left (\dfrac{p+1}{p-1}\right )^4 L(p,s) & =\left (\dfrac{p+1}{p-1} s\right )^4-\dfrac{32p(p+1)^3}{(p-1)^4}s^2+\dfrac{32p(p+3)(p+1)^2}{(p-1)^4}s-\dfrac{64p(p+1)^2}{(p-1)^4}\nonumber \\[0.15cm]
& =\mu^4-\dfrac{32p(p+1)}{(p-1)^2}\mu^2+\dfrac{32p(p+1)(p+3)}{(p-1)^3}\mu-\dfrac{64p(p+1)^2}{(p-1)^4}\\
& =: H(p,\mu). \nonumber
\end{align}Thus (\ref{eq:3.1}) is equivalent to
\begin{equation}\label{eq:3.7}
H(p,\mu)<0.
\end{equation}\vskip .08in

\noindent (ii) $\Longleftrightarrow$ (iii).
We recall from Theorem 1 in \cite{Karageorgis} that the radial entire solutions to the biharmonic equation $\Delta^2 u=u^p$ in $\mathbb{R}^N$ are unstable if and only if
\begin{equation}\label{eq:3.8}
\dfrac{N^2(N-4)^2}{16} <p K(\lambda),
\end{equation}where $\lambda =-\dfrac{4}{p-1}$ and $K(\lambda)=\lambda(\lambda-2)(\lambda+N-2)(\lambda+N-4)$. We note that the left hand side of (\ref{eq:3.8}) comes from the best constant of the Hardy-Rellich inequality (see \cite{Rellich}): If $N \ge 5$, then, for all $\varphi \in H^2(\mathbb{R}^N)$
\begin{equation*}
\int_{\mathbb{R}^N} |\Delta \varphi|^2 dx \ge \dfrac{N^2(N-4)^2}{16}\int_{\mathbb{R}^N} \dfrac{\varphi^2}{|x|^4} dx.
\end{equation*}Solving the corresponding quartic inequality, we find that
\begin{equation*}
(\ref{eq:3.8}) \;\ \mbox{holds}\; \mbox{true}\; \mbox{if}\; \mbox{and} \; \mbox{only}\; \mbox{if}\;\ 1<p<p(N),
\end{equation*}where $p(N)$ is the fourth-order Joseph-Lundgren exponent computed by Gazzola and Grunau (see \cite{Gazzola}):
\begin{equation*}
p(N)=
\begin{cases}
+\infty, & \mbox{if}\; N\le 12\\[0.1cm]
\dfrac{N+2-\sqrt{N^2+4-4\sqrt{N^2+H_N}}}{N-6-\sqrt{N^2+4-4\sqrt{N^2+H_N}}}, & \mbox{if}\; N \ge 13,
\end{cases}
\end{equation*}and $H_N=\dfrac{N^2(N-4)^2}{16}$.\vskip .05in

If we denote $N:=2+2\mu$ in (\ref{eq:3.8}), a direct calculation shows
\begin{equation*}
\mu^4-2\mu^2+1 <\dfrac{32p(p+1)}{(p-1)^2}\left [\mu^2-\dfrac{p+3}{p-1}\mu+\dfrac{2(p+1)}{(p-1)^2} \right ].
\end{equation*}Thus, it implies that (\ref{eq:3.8}) is equivalent to
\begin{equation}\label{eq:3.9}
H_0(p,\mu):=(\mu^2-1)^2-\dfrac{32p(p+1)}{(p-1)^2}\mu^2+\dfrac{32p(p+1)(p+3)}{(p-1)^3}\mu-\dfrac{64p(p+1)^2}{(p-1)^4}<0.
\end{equation}We find that
\begin{equation*}
H_0(p,\mu)=H(p,\mu)-2\mu^2+1.
\end{equation*}Combining the above identities with (\ref{eq:3.6})-(\ref{eq:3.9}), we obtain that (\ref{eq:3.7}) is equivalent to
\begin{equation*}
\dfrac{\ell^2(\ell-4)^2}{16}+\dfrac{(\ell-2)^2}{2}-1<pK(\lambda),
\end{equation*}where $\ell =2+2\mu$.

Now, we denote $H^*_{\ell}:=\dfrac{\ell^2(\ell-4)^2}{16}+\dfrac{(\ell-2)^2}{2}-1$, then we conclude that (\ref{eq:3.7}) holds true if and only if $1<p<p_*(\ell)$, where
\begin{equation*}p_*(\ell)=
\begin{cases}
+\infty,& 5 \le \ell \le \overline{\ell},\\[0.1cm]
\dfrac{\ell+2-\sqrt{\ell^2+4-4\sqrt{\ell^2+H^*_{\ell}}}}{\ell-6-\sqrt{\ell^2+4-4\sqrt{\ell^2+H^*_{\ell}}}}, & \ell >\overline{\ell},
\end{cases}
\end{equation*}and $\overline{\ell}$ is a unique root in the interval $(12,13)$ such that
$8(\ell-2)(\ell-4)=H^*_{\ell}$.
\end{proof}\vskip .1in

\begin{remark}\label{eq:r3.2}
\begin{itemize}
\item [\rm (i).] By Lemma 2.2 in \cite{Hajlaoui}, we find that, for any $p>1$,
$L(p,2t_0^+)<0$. Moreover, $L$ defined by (\ref{eq:3.1}) has a unique root $s_0$ in the interval $[2t_0^+,+\infty)$ such that $L(p,s)<0$ for $s \in [2t_0^+,s_0)$. Here
\begin{equation*}
t_0^+:=\sqrt{\dfrac{2p}{p+1}}+\sqrt{\dfrac{2p}{p+1}-\sqrt{\dfrac{2p}{p+1}}}.
\end{equation*}
\item [\rm (ii).] Let $N \ge 5$, $p>1$, $\alpha > 0$ and $(u,v)$ be a classical positive stable solution of (\ref{eq:1.1}). Adopting the similar proof as Lemma 4 in \cite{Cowan}, we find
that Lemma \ref{eq:l3.1} holds true for $s \in (2t_0^-,2t_0^+)$, where
\begin{equation*}
t_0^-:=\sqrt{\dfrac{2p}{p+1}}-\sqrt{\dfrac{2p}{p+1}-\sqrt{\dfrac{2p}{p+1}}}.
\end{equation*}Obviously, $t_0^-$ is strictly decreasing in $p$, $\lim\limits_{p \to 1} t_0^- =1$, $\lim\limits_{p \to \infty} t_0^-=\sqrt{2}-\sqrt{2-\sqrt{2}}$. Thus, for any $p>1$, $2t_0^-<2$.
\item [\rm (iii).] If $N \ge 5$, $p>1$, $\alpha > 0$ and $(u,v)$ is a classical nonnegative stable solution of (\ref{eq:1.1}), then Lemma \ref{eq:l3.1} holds true for $s \in [1,s_0)$.
\end{itemize}
\end{remark}\vskip .1in

\begin{lemma}\label{eq:l3.2}
Let $N \ge 5$, $p>1$, $\alpha > 0$ and $(u,v)$ be a classical stable solution of (1.1) with $u>0$. Assume that $1\le \beta <\dfrac{N}{N-2}$, $2t_0^- < \tau \le s$, where $s$ is defined in Lemma \ref{eq:l3.1}. Then there exist $k \in \mathbb{N}^+$ and $C<+\infty$ such that
\begin{equation*}
\left (\int_{B_{R_k}} v^{\beta \tau} dx\right )^{\frac{1}{\beta}}
\le CR^{N(\frac{1}{\beta}-1)} \int_{B_{R_{k+3}}} v^{\tau}dx,
\end{equation*}for all $R>0$.
\end{lemma}

\begin{proof}
Denote $w:=v^{\tau}$. A simple calculation yields
\begin{align*}
|\Delta w|\le \tau(\tau-1)v^{\tau-2}|\nabla v|^2+\tau(1+|x|^2)^{\frac{\alpha}{2}}v^{\tau-1}u^p.
\end{align*}Then, it implies from Lemma \ref{eq:l2.6} that
\begin{align}\label{eq:3.10}
\Big (\int_{B_{R_k}} v^{\beta \tau}dx \Big )^{\frac{1}{\beta}} \le & C R^{2+N(\frac{1}{\beta}-1)} \int_{B_{R_{k+1}}}\left [ v^{\tau-2}|\nabla v|^2+ (1+|x|^2)^{\frac{\alpha}{2}}u^pv^{\tau-1} \right ]dx \nonumber \\[0.1cm]
& +C R^{N(\frac{1}{\beta}-1)} \int_{B_{R_{k+1}}} v^{\tau} dx.
\end{align}

Next, we estimate the term $\displaystyle \int_{B_{R_{k+1}}} v^{\tau-2}|\nabla v|^2 dx$ in (\ref{eq:3.10}). We take a cut-off function $\phi \in C_0^2 (B_{R_{k+2}})$ such that  $\phi \equiv 1$ in $B_{R_{k+1}}$ and $|\nabla \phi| \le \dfrac{C}{R}$. Multiply $-\Delta v=(1+|x|^2)^{\frac{\alpha}{2}}u^p$ by $v^{\tau-1}\phi^2$ and integrate by parts to get
\begin{align}\label{eq:3.11}
\int_{\mathbb{R}^N} (1 & +|x|^2)^{\frac{\alpha}{2}} u^p v^{\tau-1}\phi^2 dx =\int_{\mathbb{R}^N} -\Delta v (v^{\tau-1}\phi^2)dx \nonumber \\[0.1cm]
= & (\tau-1)\int_{\mathbb{R}^N} v^{\tau-2}\phi^2|\nabla v|^2dx +2\int_{\mathbb{R}^N} v^{\tau-1}\phi \nabla v \cdot \nabla \phi dx.
\end{align}We use Young's inequality to yield
\begin{align*}
\int_{\mathbb{R}^N} v^{\tau-1}\phi \nabla v\cdot \nabla \phi dx \le \dfrac{1}{8} \int_{\mathbb{R}^N} v^{\tau-2}\phi^2 |\nabla v|^2 dx +C \int_{\mathbb{R}^N} v^{\tau}|\nabla \phi|^2 dx.
\end{align*}Substituting the above inequality into (\ref{eq:3.11}), it implies from the choice of the function $\phi$ that
\begin{align}\label{eq:3.12}
\int_{B_{R_{k+1}}} v^{\tau-2}|\nabla v|^2 dx \le C\int_{\mathbb{R}^N} (1+|x|^2)^{\frac{\alpha}{2}} u^p v^{\tau-1} \phi^2 dx +C\int_{\mathbb{R}^N} v^{\tau}|\nabla \phi|^2 dx \nonumber \\[0.1cm]
\le C\int_{B_{R_{k+2}}} (1+|x|^2)^{\frac{\alpha}{2}} u^pv^{\tau-1}dx +\dfrac{C}{R^2} \int_{B_{R_{k+2}}} v^{\tau} dx.
\end{align}Putting (\ref{eq:3.12}) into (\ref{eq:3.10}) gives
\begin{align}\label{eq:3.13}
\left (\int_{B_{R_k}} v^{\beta \tau} dx\right )^{\frac{1}{\beta}} \le & C R^{2+N(\frac{1}{\beta}-1)} \int_{B_{R_{k+2}}} (1+|x|^2)^{\frac{\alpha}{2}} u^pv^{\tau-1} dx \nonumber \\[0.1cm]
& +CR^{N(\frac{1}{\beta}-1)} \int_{B_{R_{k+2}}} v^{\tau} dx.
\end{align}Applying Lemma \ref{eq:l3.1} and Remark \ref{eq:r3.2} (ii) to the first term on the right hand side of (\ref{eq:3.13}), we
conclude that
\begin{equation*}
\left (\int_{B_{R_k}} v^{\beta \tau} dx\right )^{\frac{1}{\beta}}
\le CR^{N(\frac{1}{\beta}-1)} \int_{B_{R_{k+3}}} v^{\tau}dx,
\end{equation*}for all $R>0$.
\end{proof}

Now, we can follow exactly the iteration process in \cite[Corollary 2]{Cowan} to get

\begin{proposition}\label{eq:p3.1}
Suppose that $N \ge 5$, $p>1$, $\alpha>0$ and $(u,v)$ is a classical stable solution of (1.1) with $u>0$. Let $2t_0^-<\tau<2$ and $\tau \le \beta<\dfrac{N}{N-2}s$, then there is a $C<+\infty$ and integer $n \ge 1$ such that for all $R>1$
\begin{equation*}
\left (\int_{B_R}v^{\beta} dx\right )^{\frac{1}{\beta}} \le CR^{N(\frac{1}{\beta}-\frac{1}{\tau})} \left (\int_{B_{R_{3^n}}} v^{\tau}dx \right )^{\frac{1}{\tau}}.
\end{equation*}Here $s$ is defined in Lemma \ref{eq:l3.1}.
\end{proposition}

\noindent {\it Proof of Theorem \ref{eq:t1.1}.} Suppose that $(u,v)$ is a classical positive stable solution of (\ref{eq:1.1}). By Lemma \ref{eq:l3.1} and Remark \ref{eq:r3.1}-\ref{eq:r3.2}, we take $2t_0^-<\tau<2$ and $\tau \le \beta<\dfrac{N}{N-2}s$, then, combining Proposition \ref{eq:p3.1}, we obtain
\begin{equation}\label{eq:3.14}
\left (\int_{B_R} v^{\beta} dx \right )^{\frac{1}{\beta}} \le CR^{N\left (\frac{1}{\beta}-\frac{1}{\tau} \right )} \left (\int_{B_{R_{3^n}}} v^{\tau} dx \right )^{\frac{1}{\tau}}.
\end{equation}From Lemma \ref{eq:l2.4}, we apply H\"{o}lder's inequality to get
\begin{align*}
\int_{B_R} v^{\tau} dx \le & \left (\int_{B_R} (1+|x|^2)^{\frac{\alpha}{2}}v^2 dx \right )^{\frac{\tau}{2}} \left (\int_{B_R} (1+|x|^2)^{-\frac{\alpha\tau}{2(2-\tau)}}dx  \right )^{\frac{2-\tau}{2}} \\[0.1cm]
\le & CR^{\left (N-4-\alpha-\frac{8+4\alpha}{p-1} \right )\frac{\tau}{2}}R^{\left (N-\frac{\alpha \tau}{2-\tau} \right ) \frac{2-\tau}{2}}\\[0.1cm]
= & CR^{N-\frac{(2+\alpha)(p+1)}{p-1}\tau}.
\end{align*}Combine the above inequality with (\ref{eq:3.14}) to yield
\begin{equation*}
\left (\int_{B_R} v^{\beta} dx \right )^{\frac{1}{\beta}} \le CR^{N\left (\frac{1}{\beta}-\frac{1}{\tau}\right )+\frac{1}{\tau}\left (N-\frac{(2+\alpha)(p+1)}{p-1}\tau \right )},
\end{equation*}for all $R>0$. We easily find the fact that
\begin{equation*}
N\left (\frac{1}{\beta}-\dfrac{1}{\tau}\right )+\frac{1}{\tau}\left (N-\dfrac{(2+\alpha)(p+1)}{p-1}\tau \right )<0
\end{equation*}is equivalent to
\begin{equation}\label{eq:3.15}
N < \dfrac{(2+\alpha)(p+1)}{p-1} \beta  < \dfrac{(2+\alpha)(p+1)}{p-1} \dfrac{N}{N-2}s.
\end{equation}Since $\mu =\dfrac{p+1}{p-1}s$ and $\ell=2+2\mu$, it implies from Remark \ref{eq:l3.1} that the inequality (\ref{eq:3.15}) is equivalent to one of the following two conditions:
\begin{itemize}
\item [\rm (i).] $\ell \ge 5$, $p \in (1,p_*(\ell))$ and
\begin{equation*}
 N<\ell + \dfrac{\alpha(\ell-2)}{2}.
\end{equation*}
\item [\rm (ii).] For any $p>1$,
\begin{equation*}
N<2+(2+\alpha) \mu,
\end{equation*}where $\mu$ is the largest root of the polynomial
\begin{equation*}
H(p,\mu)=\mu^4-\dfrac{32p(p+1)}{(p-1)^2}\mu^2+\dfrac{32p(p+1)(p+3)}{(p-1)^3}\mu-\dfrac{64p(p+1)^2}{(p-1)^4}.
\end{equation*}
\end{itemize}
From the condition of Theorem \ref{eq:t1.1} and Remark \ref{eq:r1.1} (iii), it is easy to see that (\ref{eq:3.15}) holds true.
Therefore, we conclude that $\|v\|_{L^{\beta}(\mathbb{R}^N)} =0$ as $R \to +\infty$, i.e., $v \equiv 0$ in $\mathbb{R}^N$.
This is a contradiction. Thus we get the desired result. \hspace*{161pt} $\square$\vskip .1in

Adopting the similar proof as Theorem \ref{eq:t1.1}, we obtain the following result.

\begin{corollary}
Let $N \ge 5$, $\alpha>0$ and $p>1$ satisfy one of the following conditions:
\begin{itemize}
\item [\rm (i).] $N <\ell +\dfrac{\alpha (\ell-2)}{2}$ and $p \in (1,p_*(\ell))$. or
\item [\rm (ii).] For any $p>1$, $N<2+(2+\alpha)\mu$ and $\mu$ is the largest root of the polynomial $H(p,\mu)$ in (\ref{eq:1.5}).
\end{itemize}
Then there does no exist classical positive
stable solution of (\ref{eq:2.9}) in $\Sigma =\mathbb{R}_+^N$.
\end{corollary}\vskip .2in

\section{Proof of Theorem \ref{eq:t1.2}}

From Remark \ref{eq:l3.2} (i) and (iii), we adopt the same proof as in Lemma 3.1 and Lemma 3.2 to obtain the following results.

\begin{lemma}\label{eq:l4.1}
Let $N \ge 5$, $p>1$, $\alpha>0$ and $(u,v)$ be a classical nonnegative stable solution of (1.1). Assume that $1\le \beta <\dfrac{N}{N-2}$ and $s$ is defined in Lemma \ref{eq:l3.1}. Then there exist integer $k \in \mathbb{N}^+$ and $C<+\infty$ such that
\begin{equation*}
\left (\int_{B_{R_k}} v^{\beta s} \right )^{\frac{1}{\beta}}
\le CR^{N(\frac{1}{\beta}-1)} \int_{B_{R_{k+3}}} v^sdx,
\end{equation*}for all $R>0$.
\end{lemma}

\begin{proposition}\label{eq:p4.1}
Suppose that $N \ge 5$, $p>1$, $\alpha>0$ and $(u,v)$ is a classical nonnegative stable solution of (1.1). Then, for any $2 \le \beta <\dfrac{N}{N-2} s$, there is some $C<+\infty$ and integer $n \ge 1$ such that for all $R>1$
\begin{equation*}
\left (\int_{B_R}v^{\beta} \right )^{\frac{1}{\beta}} \le CR^{N(\frac{1}{\beta}-\frac{1}{2})} \left (\int_{B_{R_{3^n}}} v^2 \right )^{\frac{1}{2}}.
\end{equation*}Here $s$ is defined in Lemma \ref{eq:l3.1}.
\end{proposition}

\noindent {\it Proof of Theorem \ref{eq:t1.2}.}
(1). Suppose the condition (i) holds. Let $2 \le \beta <\dfrac{N}{N-2}s$, then, combining Proposition \ref{eq:p4.1} with Lemma \ref{eq:l2.4}, we get
\begin{align*}
\left (\int_{B_R} v^{\beta} dx \right )^{\frac{1}{\beta}} & \le CR^{\frac{N}{2}(\frac{2}{\beta}-1)} \left (\int_{B_{R_{3^n}}} v^2 dx \right )^{\frac{1}{2}}\\[0.1cm]
& \le CR^{\frac{N}{2}(\frac{2}{\beta}-1)}\left (\int_{B_{R_{3^n}}}(1+|x|^2)^{\frac{\alpha}{2}} v^2 dx \right )^{\frac{1}{2}}\\[0.1cm]
& \le CR^{\frac{N}{2}(\frac{2}{\beta}-1)+\frac{N}{2}-2-\frac{\alpha}{2}-\frac{4+2\alpha}{p-1}},
\end{align*}for all $R>0$. We obtain that
\begin{equation*}
\dfrac{N}{2}\left (\dfrac{2}{\beta}-1\right )+\frac{N}{2}-2-\dfrac{\alpha}{2}-\dfrac{4+2\alpha}{p-1}<0
\end{equation*}is equivalent to
\begin{equation}\label{eq:4.1}
N < \left [\dfrac{2(p+1)}{p-1}+\dfrac{\alpha(p+3)}{2(p-1)} \right ]\dfrac{N}{N-2}s.
\end{equation}Since $\mu =\dfrac{p+1}{p-1}s$ and $\ell=2+2\mu$, it implies from Remark \ref{eq:r3.1} that the inequality (\ref{eq:4.1}) is equivalent to
\begin{equation*}
 N<\ell + \dfrac{\alpha(\ell-2)(p+3)}{4(p+1)},
\end{equation*}where $\ell \ge 5$ and $p \in (1,p_*(\ell))$. Therefore, if $\ell \ge 5$, $p \in (1,p_*(\ell))$ and $N<\ell + \dfrac{\alpha(\ell-2)(p+3)}{4(p+1)}$, we obtain that $\|v\|_{L^{\beta}(\mathbb{R}^N)} =0$ as $R \to +\infty$, i.e., $v \equiv 0$ in $\mathbb{R}^N$.
Thus we get the desired result. \vskip .08in

(2). Take $\mu =\dfrac{p+1}{p-2}s$. We find that the condition (ii) is equivalent to (\ref{eq:4.1}). Hence, we get the desired result by adopting the same proof as the above.
\hspace*{90pt} $\square$\vskip .2in

\noindent {\bf Acknowledgments} \vskip .1in

The work was partially supported by NSFC of China (No. 11201248), K.C. Wong Fund of Ningbo University and
Ningbo Natural Science Foundation (No. 2014A610027).

\end{document}